# Multifidelity Space Mission Planning and Infrastructure Design Framework for Space Resource Logistics[1]

Hao Chen,[2] and Tristan Sarton du Jonchay,[3]
*Georgia Institute of Technology, Atlanta, GA, 30332*

Linyi Hou,[4]
*University of Illinois at Urbana-Champaign, Urbana, IL, 61801*

and
Koki Ho[5]
*Georgia Institute of Technology, Atlanta, GA, 30332*

**To build a sustainable space transportation system for human space exploration, the design and deployment of space infrastructure, such as in-situ resource utilization plants, in-orbit propellant depots, and on-orbit servicing platforms, are critical. The design analysis and trade studies for these space infrastructure systems require the consideration of not only the design of the infrastructure elements themselves, but also their supporting systems (e.g., storage, power) and logistics transportation while exploring various architecture options (e.g., location, technology). This paper proposes a system-level space infrastructure and logistics design optimization framework to perform architecture trade studies. A new space-infrastructure logistics optimization problem formulation is proposed that considers the internal interactions of infrastructure subsystems and their external synergistic effects with space logistics simultaneously. Because the full-size version of this proposed problem formulation can be computationally prohibitive, a new multifidelity optimization formulation is developed by varying the granularity of the commodity-type definition over the space logistics network; this multifidelity formulation can find an approximate solution to the full-size problem computationally efficiently with little sacrifice in the solution quality. The**

---

[2] Ph.D. Student, Aerospace Engineering, AIAA Student Member.
[3] Ph.D. Student, Aerospace Engineering, AIAA Student Member.
[4] Undergraduate Student, Aerospace Engineering, AIAA Student Member.
[5] Assistant Professor, Aerospace Engineering, AIAA Senior Member.

**proposed problem formulation and method are applied to the design of in situ resource utilization systems in a multimission lunar exploration campaign to demonstrate their values.**

## Nomenclature

| | | |
|---|---|---|
| $A$ | = | coefficient matrix for linear programming |
| $\mathcal{A}$ | = | set of arcs |
| $b$ | = | constraint vector for linear programming |
| $\mathcal{C}_c$ | = | continuous commodity set |
| $\mathcal{C}_d$ | = | discrete commodity set |
| $C_v$ | = | spacecraft payload capacity |
| $c$ | = | cost coefficient |
| $d$ | = | mission demand, kg |
| $F$ | = | commodity transformation matrix |
| $G$ | = | aggregation matrix |
| $g_0$ | = | standard gravity, m/s$^2$ |
| $H$ | = | concurrency constraint matrix |
| $I_{sp}$ | = | specific impulse, s |
| $\mathcal{J}$ | = | optimization objective |
| $\mathcal{N}$ | = | set of nodes |
| $P_I$ | = | infrastructure power demand, kW |
| $P_v$ | = | spacecraft propellant capacity, kg |
| $P_0$ | = | power system output power, kW |
| $Q_I$ | = | infrastructure operating length per solar day, h |
| $Q_p$ | = | power system working time per solar day, h |
| $Q_u$ | = | commodity packing index subsets |
| $R$ | = | number of different types of commodities |
| $S_k$ | = | commodity index partition subsets |
| $S_v$ | = | spacecraft structure mass, kg |

| $\mathcal{T}$ | = | set of time steps |
| --- | --- | --- |
| $\mathcal{V}$ | = | set of spacecraft |
| $\mathcal{W}$ | = | set of time windows |
| $\mathcal{X}$ | = | variable vector for linear programming |
| $x$ | = | commodity variable |
| $\alpha$ | = | reactor productivity, kg/h/kg plant |
| $\beta$ | = | consumption rate, kg/h/kg plant |
| $\varepsilon$ | = | energy storage efficiency |
| $\gamma$ | = | the specific mass of energy storage systems |
| $\Delta t$ | = | time of flight, day |
| $\Delta V$ | = | change of velocity, km/s |
| $\tilde{\zeta}$ | = | commodity packing index set |
| $\phi$ | = | propellant mass fraction |

*Subscripts*

| $i$ | = | node index |
| --- | --- | --- |
| $j$ | = | node index |
| $t$ | = | time step index |
| $v$ | = | spacecraft index |

## I. Introduction

AS interest grows in space exploration and space economy development, the design and deployment of space infrastructure systems become critical to support space resource utilization, on-orbit servicing, and interplanetary space transportation. Past space infrastructure design literature has analyzed the system performance of in-situ resource utilization (ISRU) systems, propellant depots, on-orbit servicing platforms, etc. For example, the technical and economic feasibility of commercial propellant production by ISRU systems has been examined and demonstrated by industry, government, and academic experts [1]. Multiple studies have focused on the chemical processes of ISRU reactor and system productivity, such as the hydrogen reduction reaction testbeds by NASA [2] and Lockheed Martin [3], the integrated carbothermic regolith reduction system by Orbitec Inc. and Kennedy Space Center [4], the integrated molten regolith electrolysis (MRE) reactor modeling method by Schreiner [5], and the Mars oxygen ISRU experiment by Meyen [6]. Besides ISRU, on-orbit servicing technologies have also been developed in recent years [7,8], its commercial potentials and operations have been analyzed in the literature [9-11]. However, all these referenced studies mainly analyzed the feasibility and performance of the space infrastructure elements, and did not take into account the complex logistics to deploy and support those infrastructure systems.

On the other hand, multiple studies focused on space transportation analysis and considered space infrastructure design, such as ISRU systems, together with space transportation system design. United Launch Alliance proposed the Cislunar-1000 project to build a sustainable space economy by taking advantage of lunar water ISRU plants to produce oxygen and hydrogen [12]. A series of network-based space logistics optimization methods were proposed by Ishimatsu [13], Ho [14], and Chen [15] to solve mission planning, space infrastructure design, and spacecraft design problems concurrently. Their results showed the long-term benefits of ISRU systems and propellant depots to space exploration campaigns. However, in these traditional space logistics optimization methods, referred to as the prefixed space infrastructure optimization formulation in this paper, the space infrastructure was considered as a black box, and the subsystem interactions and mass ratios were determined before taking into account space logistics optimization [13-17]. They ignored the interaction between infrastructure subsystems and space transportation mission planning.

Due to the inadequate trade studies between space infrastructure design and space transportation planning, conventional prefixed space mission planning and infrastructure design have only explored a limited trade space. For example, considering the ISRU system as a black box model would miss the tradeoff between the frequency of logistics missions and its impact on ISRU storage system size. Namely, frequent transportation missions require smaller storage subsystems but higher operation cost and complexity; whereas infrequent transportation missions require larger storage subsystems, which can also lead to higher infrastructure deployment cost. The consideration of this tradeoff requires the modeling of ISRU infrastructure subsystems and its coupling with the logistics planning. Furthermore, prefixing the space infrastructure design for architecture design can also miss the synergistic effect of space infrastructure technologies and the combination of subsystems to achieve a hybrid system design, particularly when different infrastructure technologies have common supporting subsystems. An example is an ISRU plant based on the reverse water gas shift reaction (RWGS) and Sabatier reaction (SR). The RWGS process can be used concurrently with SR to produce sufficient oxygen so that the generated $O_2$ and $CH_4$ can be used together as propellants optimally. (The oxygen/methane bipropellant has been widely considered as a propellant option to support future robotic and human missions in conjunction with ISRU systems [18,19].) Because of the similar reactants and reaction environment, the RWGS process and the SR process can share the same acquisition subsystem (for $CO_2$), liquefication & storage subsystem (for $O_2$), and power subsystem. Thus, the SR ISRU and the RWGS ISRU need to be designed together for optimal performance, and this design solution depends on the mission scenarios and the logistics planning (e.g., launch frequency, vehicle type/size, available resources from the ground or other sources). To resolve this challenge effectively, a general design optimization framework and its methods need to be developed to handle the synergistic effect of space infrastructure subsystems and the logistics system.

To effectively evaluate the impacts of the space infrastructure design to space missions with higher fidelity (i.e., considering both system-level and subsystem-level tradeoffs), we propose an interdisciplinary space infrastructure optimization framework and its optimization methods, leveraging network-based space logistics modeling. The proposed framework enables an integrated architecture trade study for future space infrastructure, considering the coupling between the subsystems design and corresponding logistics planning.

Our proposed framework has four technical innovations. First, we propose an interdisciplinary space infrastructure optimization formulation that considers infrastructure subsystems' internal interactions and their external synergistic effects with space logistics transportation simultaneously. This is a new problem in space logistics for high-fidelity space infrastructure trade studies. Second, since the full-size version of this proposed problem formulation can be computationally prohibitive for large-scale space infrastructure design problems, we develop a new multi-fidelity optimization formulation that can provide an approximate solution to the full-size formulation at a significantly reduced computational cost with little sacrifice in the solution quality. The idea behind this multi-fidelity formulation is to vary the granularity of the commodity type definition over the network graph; this technique is referred to as commodity packing based on its physical meaning. Third, in order to identify where and when to pack the commodities for the multi-fidelity optimization formulation, we develop a preprocessing algorithm for commodity packing. This method enables an automated implementation of the multi-fidelity formulation in the context of dynamic generalized multicommodity network flow. Fourth, we establish the relationship between the solutions of the multi-fidelity, full-size, and traditional prefixed formulations. This relationship enables us to find the approximate solution of the computationally prohibitive full-size formulation with the knowledge about the worst possible approximation error.

Our method will enable a unique tradeoff that could not be performed with traditional methods. The proposed framework can perform space infrastructure technology selection and system sizing for each subsystem considering their interactions and logistics mission planning. For example, we can consider the tradeoff between the frequency of logistics missions and its impact on ISRU storage system size. Exploring this tradeoff considering both ISRU and logistics mission design concurrently can lead to an efficient system design compared with the traditional methods considering these two separately. In addition, our method can consider the optimal design of a hybrid infrastructure system with multiple technologies. For example, for ISRU infrastructure systems that have common reactants, reaction environments, or final products like the aforementioned SR/RWGS example, the proposed framework can combine multiple technologies into an optimally integrated ISRU architecture with shared subsystems. By enabling these new capabilities in space logistics optimization, the developed framework provides an important step forward in integrated space infrastructure design and trade studies for future large-scale human space exploration. The proposed optimization framework can also be used as an evaluation tool to analyze the long-term performance of spacecraft and space architectures.

The remainder of this paper is organized as follows. Section II first introduces the traditional prefixed optimization formulation for space infrastructure design, where space infrastructure is considered as a black box. Then, Sec. III discusses the full-size version of the proposed space infrastructure optimization problem formulation, taking into account space infrastructure subsystems tradeoffs together with space mission planning concurrently. In Sec. IV, we propose a multi-fidelity optimization formulation and its methods to resolve the computational challenge inherent in the full-size formulation. Section V demonstrates the proposed optimization formulations through a multi-mission human lunar exploration campaign case study. Finally, Sec. VI summarizes the conclusion of this paper and discusses future work.

## II. Traditional Method: Prefixed Space Infrastructure Optimization Formulation

The network-based space logistics optimization formulation considers space missions as commodity network flow problems [13-15], where nodes represent planets or orbits and arcs represent trajectories. Vehicles, payloads, infrastructure, and crewmembers are all considered as commodities. An example of the Earth-Moon-Mars transportation network model is shown in Fig. 1. The inputs of this infrastructure optimization formulation are space mission demands and corresponding available infrastructure systems to be implemented (i.e., mainly representing ISRU systems and their supporting structures in this paper). Based on the mission demands and time window constraints, this formulation outputs selected infrastructure systems to be deployed, including system sizing, plant deployment strategy, system operating mechanisms, and further resource logistics processes if mission demands occur at a location different from the infrastructure deployment spot.

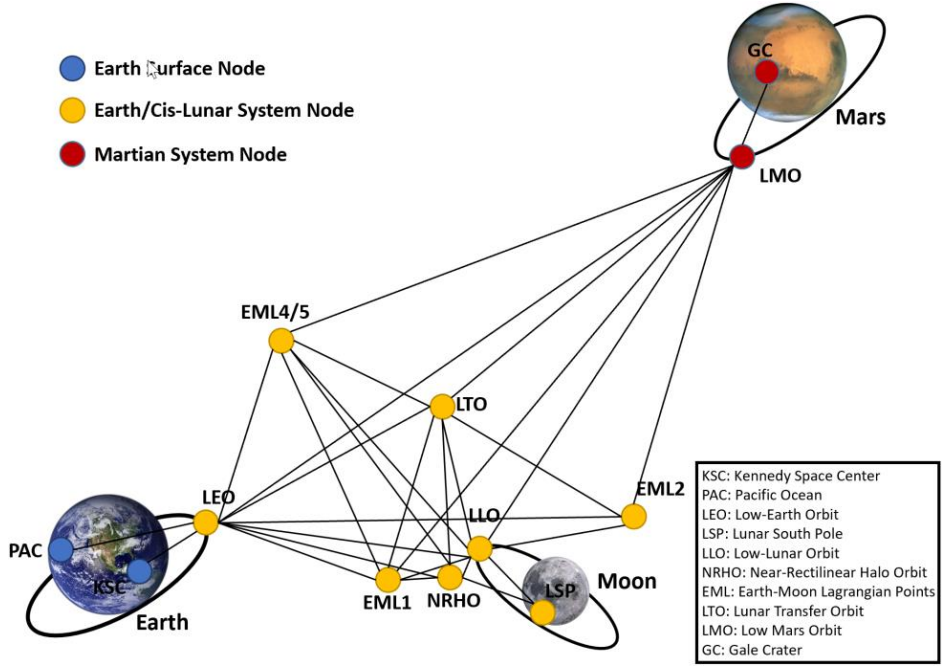

**Fig. 1 An example of the Earth-Moon-Mars transportation network model. [17]**

In this space infrastructure optimization problem, space logistics mission planning is the main goal for optimization. Fig. 1 shows an example of the Earth-Moon-Mars transportation network model. Space logistics optimization includes space transportation scheduling and space infrastructure deployment strategy optimization. The space infrastructure subsystem interactions are determined in advance before space logistics optimization. The optimizer of this formulation only finds the optimal total mass of the space infrastructure, where the mass ratios between subsystems are fixed.

Let us first define a time-expanded network graph by a set of arcs, $\mathcal{A} = \{\mathcal{N}, \mathcal{T}, \mathcal{V}\}$, where $\mathcal{N}$ is a set of nodes (index: $i, j$), and $\mathcal{T}$ is a set of time steps (index: $t$). We also need a set of available spacecraft $\mathcal{V}$ (index: $v$) as transportation vehicles during space missions. There are two types of arcs in the network: 1) transportation arcs to connect different nodes at different time steps representing spaceflights in space transportation; 2) holdover arcs to connect the same nodes at different time steps representing operation activities after infrastructure deployment. Then, we define a commodity flow variable vector $\boldsymbol{x}_{vijt}$, representing the commodity flow from node $i$ to node $j$ at time $t$ using spacecraft $v$. Note that this $\boldsymbol{x}_{vijt}$ represents the mass right after departing node $i$, and thus it is often mentioned as an outflow in the literature. Each element of the commodity flow variable vector $\boldsymbol{x}_{vijt}$ corresponds to one type of commodity, and it can be either continuous or discrete (i.e., integer) depending on the corresponding commodity; the former commodity set (i.e., continuous commodity set) is defined as $C_c$, and the latter commodity set (i.e., discrete commodity set) is defined as $C_d$. For example, the number of spacecraft and crew members are integers while the mass of propellant and payload are continuous. Define a demand parameter $\boldsymbol{d}_{it}$, which is determined by mission scenarios. Mission demands are negative and mission supplies are positive. We also need to define a cost coefficient, $\boldsymbol{c}_{vijt}$, for each commodity to measure the space mission cost. If there are $R$ types of commodities in the space mission, $\boldsymbol{x}_{vijt}$, $\boldsymbol{d}_{it}$, and $\boldsymbol{c}_{vijt}$ are all $R \times 1$ vectors.

Aside from the notations defined above, we also need to define the following parameters:

1) $\Delta t_{ij}$ = Time of flight along arc $i$ to $j$.

2) $F_{vij}$ = Commodity transformation matrix.

3) $H_{vij}$ = Concurrency constraint matrix.

4) $W_{ij}$ = Time windows of spaceflight along arc $i$ to $j$.

Then, the formulation of the prefixed space infrastructure optimization formulation can be written as follows.

Minimize:

$$\mathcal{J} = \sum_{(v,i,j,t)\in\mathcal{A}} \boldsymbol{c}_{vijt}{}^T \boldsymbol{x}_{vijt} \tag{1a}$$

Subject to:

$$\sum_{(v,j):(v,i,j,t)\in\mathcal{A}} \boldsymbol{x}_{vijt} - \sum_{(v,j):(v,j,i,t)\in\mathcal{A}} F_{vji}\boldsymbol{x}_{vji(t-\Delta t_{ji})} \leq \boldsymbol{d}_{it} \quad \forall i \in \mathcal{N} \;\; \forall t \in \mathcal{T} \tag{1b}$$

$$H_{vij}\boldsymbol{x}_{vijt} \leq \boldsymbol{0}_{l\times 1} \quad \forall (v,i,j,t) \in \mathcal{A} \tag{1c}$$

$$\begin{cases} \boldsymbol{x}_{vijt} \geq \boldsymbol{0}_{R\times 1} & \text{if } t \in W_{ij} \\ \boldsymbol{x}_{vijt} = \boldsymbol{0}_{R\times 1} & \text{otherwise} \end{cases} \quad \forall (v,i,j,t) \in \mathcal{A} \tag{1d}$$

$$\boldsymbol{x}_{vijt} = \begin{bmatrix} x_1 \\ x_2 \\ \vdots \\ x_R \end{bmatrix}_{vijt}, \quad \begin{matrix} x_n \in \mathbb{R}_{\geq 0} \; \forall n \in C_c \\ x_n \in \mathbb{Z}_{\geq 0} \; \forall n \in C_d \end{matrix} \quad \forall (v,i,j,t) \in \mathcal{A}$$

### A. Objective Function
Equation (1a) is the objective function that minimizes the total mission cost throughout the whole space campaign. Different types of mission objectives can be implemented depending on the mission performance metric.

### B. Mass Balance Constraint

Equation (1b) is the mass balance constraint that guarantees mission demands are always satisfied at node $i$. This constraint contains an inequality rather than equality to allow the possibility of dumping commodities out of the logistics system.

In this constraint, the second term $F_{vji}x_{vji(t-\Delta t_{ji})}$ represents the outcome of the commodity transformation process from node $j$ to node $i$ during spaceflights or mission operations. The transformation includes propellant burning that consumes propellant during spaceflights, crew consumptions that include food, water, and oxygen, and resource productions (e.g., propellant) by space infrastructure systems. The matrix $F_{vij}$ is the transformation matrix. After the transformation process, the commodities flow into node $i$ as commodity inflows.

To illustrate the settings of the transformation matrix $F_{vij}$, two examples are shown in the following. One example is about propellant burning and another one is about space infrastructure resource productions. First, define a commodity inflow variable as,

$$x_{vijt}^{inflow} = F_{vij}x_{vijt}$$

For propellant burning process, define the commodity flow variables as,

$$x_{vijt} = \begin{bmatrix} x^C: \text{cargo, kg} \\ x^r: \text{propellant, kg} \\ x^S: \text{spacecraft, \#} \end{bmatrix}_{vijt}$$

Then, we can express the impulsive propellant consumption as follows:

$$\begin{bmatrix} x^C \\ x^r \\ x^S \end{bmatrix}_{vijt}^{inflow} = \begin{bmatrix} 1 & 0 & 0 \\ -\phi & 1-\phi & -\phi S_v \\ 0 & 0 & 1 \end{bmatrix}_{vij} \begin{bmatrix} x^C \\ x^r \\ x^S \end{bmatrix}_{vijt} \qquad (2)$$

In Eq. (2), $S_v$ is the spacecraft structure mass; note that, $x^S$ is in the unit of the number of spacecraft and needs to be converted into the mass in kilograms. The propellant mass fraction $\phi$ is defined from the rocket equation, $\phi = 1 - \exp\left(-\frac{\Delta V}{I_{sp}g_0}\right)$, where $\Delta V$ is the change of velocity for the spaceflight, $I_{sp}$ is the specific impulse, and $g_0$ is the standard gravity.

For space infrastructure resource productions, we use lunar water ISRU as an example. The water ISRU will first extract water from lunar regolith and then electrolyze water to generate $O_2$ and $H_2$. Define the commodity flow variables as

$$x_{vijt} = \begin{bmatrix} x^{O_2}, \text{kg} \\ x^{H_2}, \text{kg} \\ x^{ISRU}, \text{kg} \end{bmatrix}_{vijt}$$

Then, we can express the ISRU production process for one hour for $O_2$ and $H_2$ as follows:

$$\begin{bmatrix} x^{O_2} \\ x^{H_2} \\ x^{ISRU} \end{bmatrix}_{vijt}^{inflow} = \begin{bmatrix} 1 & 0 & \alpha_{ISRU}^{O_2} \\ 0 & 1 & \alpha_{ISRU}^{H_2} \\ 0 & 0 & 1 \end{bmatrix}_{vij} \begin{bmatrix} x^{O_2} \\ x^{H_2} \\ x^{ISRU} \end{bmatrix}_{vijt} \quad (3)$$

In Eq. (3), there are three constraints. The first two constraints represent the ISRU production for $O_2$ and $H_2$ for one hour, where $\alpha$ is the ISRU plant productivity, representing the amount of resource generation per hour per unit mass of the ISRU plant. The last constraint means the ISRU plant system mass does not change during the production process.

## C. Concurrency Constraint

Equation (1c) is the concurrency constraint denoting commodity flow bounds. In this constraint, $l$ is the number of concurrency constraints to be considered. To illustrate the settings of the concurrency constraint matrix $H_{vij}$, two examples are shown as follows. One example is the constraints from spacecraft propellant and payload capacities. Another example is the non-negativity of commodity variables. For spacecraft propellant and payload capacities, we define the commodity flow variables as,

$$x_{vijt} = \begin{bmatrix} x^C: \text{cargo, kg} \\ x^r: \text{propellant, kg} \\ x^S: \text{spacecraft, \#} \end{bmatrix}_{vijt}$$

Then, we can express the spacecraft payload and propellant capacities as follows,

$$\begin{bmatrix} 1 & 0 & -C_v \\ 0 & 1 & -P_v \end{bmatrix}_{vij} \begin{bmatrix} x^C \\ x^r \\ x^S \end{bmatrix}_{vijt} \leq \begin{bmatrix} 0 \\ 0 \end{bmatrix} \quad (4)$$

where $C_v$ and $P_v$ are the payload and propellant capacities of spacecraft $v$, respectively. Note that, both $C_v$ and $P_v$ are spacecraft design parameters that also can be considered as design variables in the optimization, which will make the problem an integrated space mission planning and spacecraft design problem. In this scenario, the concurrency constraint has quadratic terms and the spacecraft design model may also be nonlinear. This nonlinear problem can be solved effectively using a piecewise-linear approximation method as shown in Ref. [15]. For this research, the spacecraft design is not considered as part of the optimization. Thus, the values of $C_v$ and $P_v$ are both constants in the formulation.

For the non-negativity of commodity inflow variables, we have,

$$x_{vijt}^{inflow} \geq \mathbf{0}_{R \times 1}$$

which is equivalent to,

$$-F_{vij} x_{vijt} \leq \mathbf{0}_{R \times 1}$$

In this constraint, the concurrency constraint matrix $H_{vij}$ corresponds to the negative of the transformation constraint matrix, $-F_{vij}$. It guarantees the feasibility of commodity transformations during spaceflights or surface system operation.

### D. Time Window Constraint

Equation (1d) is the time window constraint on rocket launch and spaceflight. Only when the time windows are open, are spaceflights and mission operations permitted. Otherwise, the commodity flow variable is set to be zero.

### E. Limitations of the Traditional Formulation

In this traditional space infrastructure optimization formulation, the infrastructure subsystem designs are determined in advance. Space logistics optimization only identifies the optimal total size of the space infrastructure in space missions and cannot optimize the mass ratio between subsystems. It ignores the interaction between space infrastructure subsystems and space logistics transportation planning. This formulation is not able to perform sufficient trade studies for infrastructure technology selections and identify technology gaps.

## III. Full-Size Space Infrastructure Optimization Formulation

To increase the space infrastructure design fidelity and take into account the detailed interactions between space infrastructure subsystems and space logistics transportation, this section introduces a newly developed full-size space infrastructure optimization formulation that considers all infrastructure subsystems separately throughout the space campaign.

As shown in Fig. 2, there are two main components to be optimized in the full-size space infrastructure optimization formulation. The first component is the same as the prefixed space infrastructure optimization formulation, as shown on the right side of Fig. 2. It considers space transportation mission planning, space infrastructure deployment strategy, and resource logistics after production. The second component is the space infrastructure trade studies, as shown in the left side of Fig. 2. It considers the internal tradeoffs among space infrastructure subsystems and their external interactions with space transportation to provide infrastructure subsystem sizing and technology selections. In this paper, we use ISRU systems as an example of multi-subsystem space infrastructure optimization, although the proposed method can be generally implemented in different types of space infrastructure design trade studies. There are six subsystems considered in the ISRU infrastructure model: reactor, excavator (for soil) or acquisition system (for Martian atmosphere), separator, hopper/feed/secondary subsystem, storage system, and power system. There can be multiple different reactors, excavators, etc. depending on the ISRU technologies. These subsystems are all considered as different commodities in space logistics to enable effective analysis of subsystem interactions.

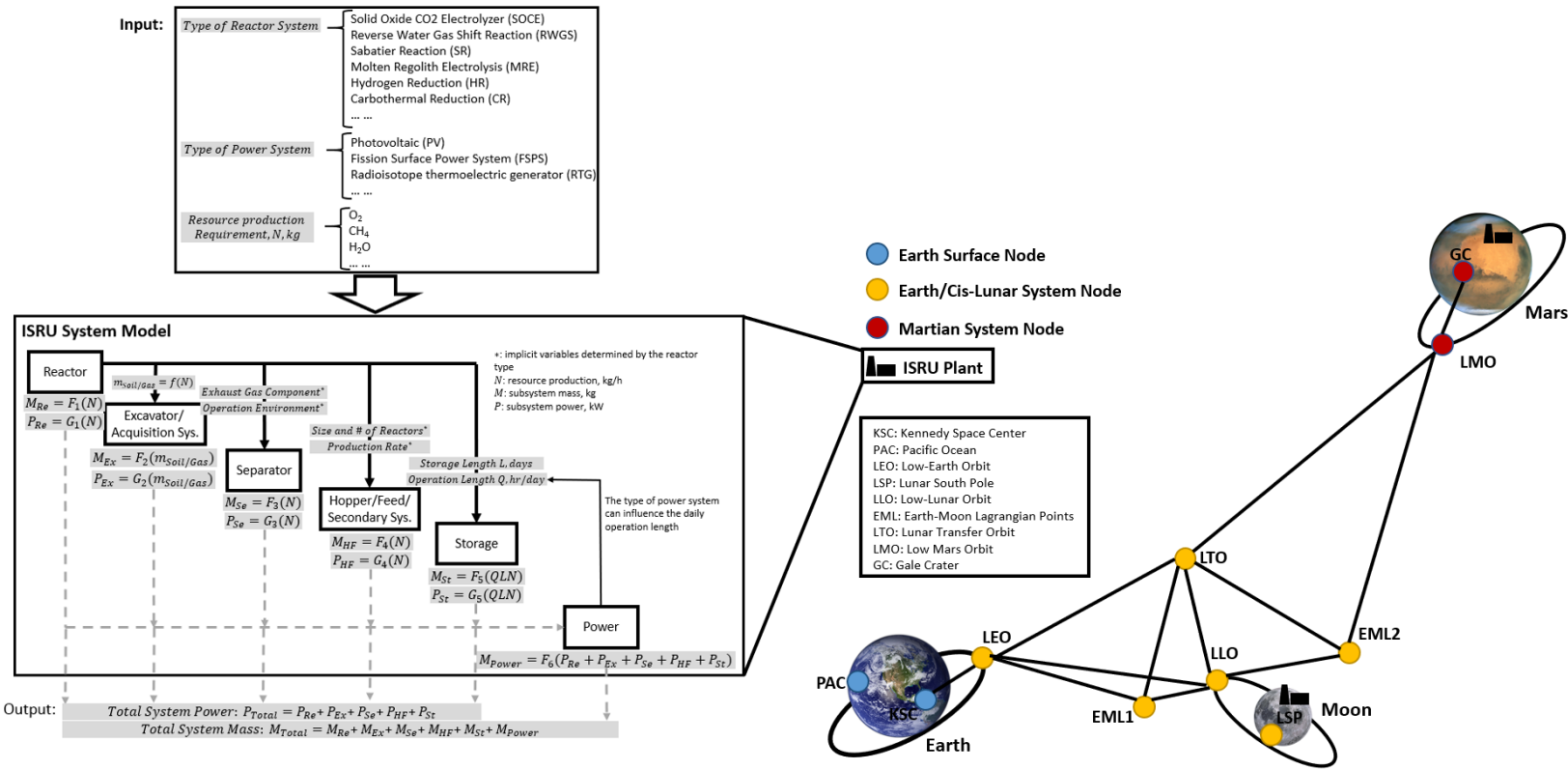

**Fig. 2 An example of the full-size space infrastructure optimization formulation.**

The formulation of the full-size space infrastructure optimization formulation is the same as the prefixed optimization formulation, as shown in Eqs. (1a)-(1d). However, the constraints are interpreted and implemented in a different way because each infrastructure subsystem is considered separately. In the following parts, we show the additional relationships to be considered to enable system-level space infrastructure trade studies together with space logistics optimization.

### A. Objective Function

The objective function (i.e., Eq. (1a)) is exactly the same as it is in the prefixed optimization formulation. The only point to note is that a higher fidelity mission performance measurement model is needed in this formulation because each subsystem is considered independently. For example, if a cost model is implemented in the objective function, then the cost model in the full-size optimization formulation should include the detailed cost information for each subsystem and each technology.

### B. Mass Balance Constraint

In the mass balance constraint (i.e., Eq. (1b)), we need to take into account the ISRU resource production process from the subsystem-level. The same lunar water ISRU example is used to illustrate the differences in the setting of the transformation matrix $F_{vij}$. There are multiple technology options to build a lunar water ISRU. In this example, we assume that the lunar water ISRU plant is mainly made up of two reactors: the soil/water extraction (SWE) reactor and the direct water electrolysis (DWE) reactor. The SWE reactor, $x^{reactor\_SWE}$, extracts water from lunar or Martian soil. The DWE reactor, $x^{reactor\_DWE}$, electrolyzes water produced by the SWE reactor to generate $O_2$ and $H_2$. We can define the commodity flow variables as,

$$\boldsymbol{x}_{vijt} = \begin{bmatrix} x^{O_2}, \text{kg} \\ x^{H_2}, \text{kg} \\ x^{H_2O}, \text{kg} \\ x^{reactor\_DWE}, \text{kg} \\ x^{reactor\_SWE}, \text{kg} \end{bmatrix}_{vijt}$$

Then, we can express the ISRU production process for one hour for $O_2$, $H_2$, and $H_2O$ as follows:

$$\begin{bmatrix} x^{O_2} \\ x^{H_2} \\ x^{H_2O} \\ x^{reactor\_DWE} \\ x^{reactor\_SWE} \end{bmatrix}^{inflow}_{vijt} = \begin{bmatrix} 1 & 0 & 0 & \alpha^{O_2}_{DWE} & 0 \\ 0 & 1 & 0 & \alpha^{H_2}_{DWE} & 0 \\ 0 & 0 & 1 & -\beta^{H_2O}_{DWE} & \alpha^{H_2O}_{SWE} \\ 0 & 0 & 0 & 1 & 0 \\ 0 & 0 & 0 & 0 & 1 \end{bmatrix}_{vij} \begin{bmatrix} x^{O_2} \\ x^{H_2} \\ x^{H_2O} \\ x^{reactor\_DWE} \\ x^{reactor\_SWE} \end{bmatrix}_{vijt} \quad (5)$$

In Eq. (5), there are five constraints in total. The first two constraints represent $O_2$ and $H_2$ generations by the DWE reactor for one hour, where $\alpha$ is the reactor productivity. The third constraint illustrates $H_2O$ consumption by the DWE reactor and production by the SWE process for one hour, where $\beta$ denotes the consumption rate. Both $\alpha$ and $\beta$ are nonnegative values. Note that because of the mass balance of chemical reactions, we have $\alpha_{DWE}^{O_2} + \alpha_{DWE}^{H_2} \leq \beta_{DWE}^{H_2O}$. The last two constraints show that the masses of the DWE reactor and the SWE reactor do not change during the resource production processes.

## C. Concurrency Constraint

In the concurrency constraint (i.e., Eq. (1c)), besides the spacecraft payload and propellant capacities considered during space transportation, the resource storage capacities for infrastructure storage systems, the power supply capacities for power generation systems, and the energy storage capacities for energy storage systems also need to be considered. Among these, the constraint format of resource storage capacities is the same as the constraints for spacecraft payload and propellant capacities, which are shown in Eq. (4).

In the following, we show two examples of the concurrency constraint in the full-size optimization formulation. One example is about space infrastructure power supply capacities and another is about energy storage capacities. For space infrastructure power supply capacities, we define the commodity flow variables as

$$\boldsymbol{x}_{vijt} = \begin{bmatrix} x^{I_1}: \text{ infrastructure system 1, kg} \\ x^{I_2}: \text{ infrastructure system 2, kg} \\ x^{I_3}: \text{ infrastructure system 3, kg} \\ x^{P}: \text{ power generation system, kg} \end{bmatrix}_{vijt}$$

Then, we can express the power supply capacity constraint for infrastructure system design as follows,

$$\left[ P_{I_1}(1 + \frac{Q_{I_1} - Q_p}{\varepsilon Q_p}) \quad P_{I_2}(1 + \frac{Q_{I_2} - Q_p}{\varepsilon Q_p}) \quad P_{I_3}(1 + \frac{Q_{I_3} - Q_p}{\varepsilon Q_p}) \quad -P_0 \right]_{vij} \begin{bmatrix} x^{I_1} \\ x^{I_2} \\ x^{I_3} \\ x^{P} \end{bmatrix}_{vijt} \leq 0 \qquad (6)$$

where $P_{I_i}, i \in \{1,2,3\}$ is the infrastructure power demand of system $i$ (in kW/kg); $Q_{I_i}, i \in \{1,2,3\}$ is the infrastructure operating length per solar day, in the unit of hours; $P_0$ is the power generation system output power per unit mass (in kW/kg); $Q_p$ is the power system working time per solar day. If the power system is a fission surface power system (FSPS) or a radioisotope thermoelectric generator (RTG), it works continuously during the space mission, which means $Q_p$ is equal to the length of a solar day. If the power system is a photovoltaic (PV) system, it only works during the daytime, which means $Q_p$ is equal to the daytime length of a solar day at the destination. If the infrastructure system operating time is longer than the power system working time per solar day, which means $Q_{I_i} > Q_p$, an energy storage system (e.g., battery or fuel cell) is necessary to support the infrastructure systems. There is an energy loss during the power storage process in battery charging/discharging. Therefore, we define an energy storage efficiency parameter, $\varepsilon$.

To identify the size of the energy storage system, a concurrency constraint for energy storage capacities is needed. Define the commodity flow variables as

$$\boldsymbol{x}_{vijt} = \begin{bmatrix} x^{I_1}: \text{infrastructure system 1, kg} \\ x^{I_2}: \text{infrastructure system 2, kg} \\ x^{I_3}: \text{infrastructure system 3, kg} \\ x^{P}: \text{power generation system, kg} \\ x^{E}: \text{energy storage system, kg} \end{bmatrix}_{vijt}$$

Then, we can express the energy storage capacity constraint as follows

$$\begin{bmatrix} -P_{I_1} & -P_{I_2} & -P_{I_3} & P_0 & -\dfrac{\gamma}{\varepsilon Q_p} \end{bmatrix}_{vij} \begin{bmatrix} x^{I_1} \\ x^{I_2} \\ x^{I_3} \\ x^{P} \\ x^{E} \end{bmatrix}_{vijt} \leq 0 \qquad (7)$$

where $\gamma$ is the specific mass of the energy storage system, in the unit of kWh/kg. It shows the ability of energy storage per unit mass.

### D. Time Window Constraint

The time window constraint (i.e., Eq. (1d)) is the same as in the prefixed optimization formulation. Typically, the time windows for different space infrastructure subsystems are the same.

### E. Relationship with the Prefixed Formulation

It is easy to show that the solution from the prefixed formulation $\mathcal{J}_{\text{prefixed}}$ is an upper bound of that from the full-size formulation $\mathcal{J}_{\text{full\_size}}$.

$$\mathcal{J}_{\text{full\_size}} \leq \mathcal{J}_{\text{prefixed}}$$

This is because the only difference between the two formulations is that the prefixed formulation fixes the mass ratios of the infrastructure subsystems, whereas the full-size formulation allows the variation of those mass ratios. Thus, the prefixed formulation has an equal or smaller feasible design space than the full-size formulation, and thus provide an equal or larger solution.

### F. Limitations of the Full-Size Formulation

The full-size space infrastructure optimization formulation considers all infrastructure subsystems as separated commodities during the entire space campaign, and this significantly increases the number of commodities in logistics optimization. Generally, as a mixed-integer linear programming formulation, which is an NP-hard problem, the computational time cost increases exponentially as the problem size increases. Studies showed that even for a short lunar exploration campaign (i.e., including three lunar missions), the concurrent optimization of space mission planning, spacecraft design, and space infrastructure design can make the network-based space logistics optimization formulation computationally prohibitive [15]. This caveat can make the full-size formulation computationally intractable for long-term space mission planning. In the next section, we will propose a new approximate optimization formulation that can achieve a significant computational cost saving with little sacrifice in the solution quality.

## IV. Multi-Fidelity Space Infrastructure Optimization Formulation

In response to the computational challenge of the full-size space infrastructure optimization formulation, we propose a new approximate optimization problem formulation. Our idea is to note the fact that the infrastructure subsystem design trade studies only exist at the destination nodes, where these subsystems are deployed; there may exist redundant commodity variables and constraints in transportation arcs that can be reduced. With this idea, we develop a mechanism to combine the infrastructure subsystem variables into fewer commodity variables during space transportation ("packing" process) and separate these packed commodities after delivery to the destination nodes ("unpacking" process). Namely, we vary the granularity of the commodity type definition over the network graph, resulting in a multi-fidelity space infrastructure optimization formulation. This formulation can significantly reduce the number of commodity variables and corresponding constraints in space logistics during space transportation and improve computational efficiency.

The multi-fidelity space infrastructure optimization leverages the theory of constraint and variable aggregations for a general mixed-integer linear programming formulation. For large and complex engineering problems, we often need to balance the accuracy of the model with the cost of computation. Constraint and variable aggregation methods have been explicitly or implicitly used in realistic problems, which are typically large and complex, to find surrogate models of the original formulations. Zipkin [20,21] performed thorough analyses on solution bounds for linear programming through constraint aggregation and variable aggregation, respectively, under certain assumptions about the problem, although their assumptions limit their methods' applicability to our problem. In the multicommodity network flow context, Evans [22,23] developed the commodity aggregation for multicommodity capacitated transportation problems to find the lower bound. More recently, Ho [24] also developed a formulation based on constraint aggregation and variable aggregation to enable an efficient way to reduce the size of the time-expanded network for the generalized multicommodity network flow.

In this section, we first discuss the general constraint and variable aggregations in linear programming. Then, we show how to perform a partial constraint and commodity aggregations, referred to as commodity packing based on its physical meaning, over particular space transportation arcs to enable a multi-fidelity optimization. We show that the solution of this multi-fidelity optimization formulation provides a lower bound of that of the full-size optimization formulation. Furthermore, a commodity packing preprocessing algorithm is also developed to enable an automatic decision on where and when to pack the commodities.

## A. Constraint Aggregation and Variable Packing

The commodity variable packing is processed in two steps: constraint aggregation and variable packing. The first step, constraint aggregation, aggregates the constraints with designated packable commodities into shared constraints through an aggregation matrix. Then, the second step, variable packing, aggregates the packable commodities into shared package commodities. The transportation, transformation, and flow bounds of these commodities are considered together through the package commodities.

### 1. Constraint Aggregation

Consider a general (full-size) linear programming formulation showed as follows.

**Formulation F1 (Full-Size)**

Minimize:
$$J = C\mathcal{X} \tag{8a}$$

Subject to:
$$A\mathcal{X} \leq b \tag{8b}$$

where

$$\mathcal{X} = \begin{bmatrix} x_1 \\ x_2 \\ \vdots \\ x_n \end{bmatrix}, C = [c_1 \quad c_2 \quad \cdots \quad c_n], A = \begin{bmatrix} a_{1,1} & a_{1,2} & \cdots & a_{1,n} \\ a_{2,1} & a_{2,2} & \cdots & a_{2,n} \\ \vdots & \vdots & \ddots & \vdots \\ a_{m,1} & a_{m,2} & \cdots & a_{m,n} \end{bmatrix}, b = \begin{bmatrix} b_1 \\ b_2 \\ \vdots \\ b_m \end{bmatrix}$$

We define an "aggregation matrix" $G$ and multiply both sides of the constraint Eq. (8b) by $G$. Then, we can obtain a new formulation as follows,

**Formulation F2 (Constraint Aggregation)**

Minimize:

$$J = C\mathcal{X} \tag{9a}$$

Subject to:

$$GA\mathcal{X} \leq Gb \tag{9b}$$

where the aggregation matrix $G$ has a size $K \times m$, where $m$ is the number of rows in the $A$ matrix and $K$ is the number of constraints after aggregation ($K \leq m$), and satisfies the following two conditions:

*Condition 1:* The aggregation matrix $G$ has exactly one nonzero entry per column, and that entry is positive.

*Condition 2:* The aggregation matrix $G$ has at least one nonzero entry per row, and those entries are all positive.

For these formulations, we show that a lower bound of the optimal objective of **F1** can be found by solving **F2** if both problems are feasible and bounded.

We first rewrite the constraint (9b) as,

$$G(A\mathcal{X} - b) \leq \mathbf{0}_{K \times 1}$$

The column indices of the positive entries in each row of the aggregation matrix $G$ define a partition of the corresponding constraints $\{1,\ldots,m\}$ into $K$ sets. Denote the partition as $\sigma = \{S_k : k = 1, \ldots, K\}$, where $S_k$ is the set of constraint indices in the $k$-th set. Define $m_k = |S_k|$, which is the number of constraint indices in the $k$-th set. The partition satisfies

$$\bigcup_{k=1}^{K} S_k = \{1, \ldots, m\} \text{ and } S_k \cap S_{k'} = \emptyset \quad \forall k \neq k'$$

Define $G = \begin{bmatrix} \boldsymbol{g}_1^T \\ \boldsymbol{g}_2^T \\ \vdots \\ \boldsymbol{g}_K^T \end{bmatrix}$, where each row of the aggregation matrix is a $1 \times m$ weighting vector, $\boldsymbol{g}_k^T$, that satisfies

$$\begin{cases} \boldsymbol{g}_k[j] > 0 & \text{if } j \in S_k \\ \boldsymbol{g}_k[j] = 0 & \text{if } j \notin S_k \end{cases} \quad \forall k \in \{1, \ldots, K\}$$

To aggregate and relax the constraints, we replace each subset of constraints $S_k$ by a single constraint through weighting vectors. As a result, we can write the $k$-th constraint after aggregation for **F1** as,

$$\boldsymbol{g}_k^T(A\mathcal{X} - b) \leq 0$$

This constraint aggregates $m_k$ number of constraints in **F1** with indices $\{j: j \in S_k\}$. Because all non-zero entries in the weighting vectors are positive, these constraints are also relaxed. By applying the weighing vectors to **F1**, we can get a relaxed formulation,

Minimize:
$$J = C\mathcal{X} \tag{10a}$$

Subject to:
$$\begin{bmatrix} g_1^T \\ g_2^T \\ \vdots \\ g_K^T \end{bmatrix} (A\mathcal{X} - b) \leq \mathbf{0}_{K \times 1} \tag{10b}$$

By solving the formulation (10a)-(10b), which is equivalent to **F2**, we can get a lower bound of **F1**'s solution.

*2. Variable Packing*

After the constraint aggregation, we can perform variable packing to further improve computational efficiency by reducing the number of variables. The purpose of this step is to find a formulation equivalent to **F2**, but with fewer variables; this step corresponds to packing the commodities. Note that, in the following discussion, we only consider the aggregation of the continuous commodity flow variables for simplicity.

Consider a variable vector as follows,

$$\mathcal{X} = \begin{bmatrix} x_1 \\ x_2 \\ \vdots \\ x_n \end{bmatrix}$$

Assume that there exists a set of index set $Q = \{Q_u: u = 1, \ldots, U\}$, where each set $Q_u$ includes the packable commodity variable indices to be packed into one package commodity $\widetilde{x_u} = \sum_{i \in Q_u} x_i$. This index set $Q$ satisfies

$$\bigcup_{u=1}^{U} Q_u \subseteq \{1, \ldots, n\} \text{ and } Q_u \cap Q_{u'} = \emptyset \quad \forall u \neq u'$$

The variable packing operation is defined as replacing the $n$ original variables $\mathcal{X}$ into $U$ new variables $\widetilde{\mathcal{X}}$ following the conversion $\widetilde{x_u} = \sum_{i \in Q_u} x_i$.

In the following, we show that we can find an equivalent formulation after performing variable packing if coefficients in **F2** satisfy the following two conditions:

*Condition 3*: For each index set $Q_u$, there exists a constant $c_u'$ such that $c_i = c_u'$ for all $i \in Q_u$;

*Condition 4*: For each index set $Q_u$, there exists a constant vector $R'_u = [r'_1, r'_2, \ldots, r'_K]^T_u$ such that $\sum_{j=1}^{m} g_{k,j} a_{j,i} = r'_k$ for all $i \in Q_u$ and for all $k \in \{1, \ldots, K\}$.

First, without loss of generality, we consider a case where the last $n - q$ variables are to be packed into one package commodity. This corresponds to the case where $Q = \{q + 1, \ldots, n\}$ and $U = 1$. Thus, the expected variable vector after packing is

$$\widetilde{\boldsymbol{X}} = \begin{bmatrix} x_1 \\ \vdots \\ x_q \\ \tilde{x} \end{bmatrix}$$

where the package commodity variable $\tilde{x} = \sum_{i=q+1}^{n} x_i$. In the objective function of **F2** (i.e., Eq. (9a)), we have,

$$C\boldsymbol{X} = [c_1 \quad c_2 \quad \cdots \quad c_n] \begin{bmatrix} x_1 \\ x_2 \\ \vdots \\ x_n \end{bmatrix} = [c_1 \quad \cdots \quad c_q \quad c_{q+1} \quad \cdots \quad c_n] \begin{bmatrix} x_1 \\ \vdots \\ x_q \\ x_{q+1} \\ \vdots \\ x_n \end{bmatrix}$$

From the first condition, we know that $c_i = c'$ for all $i \in \{q + 1, \ldots, n\}$. Therefore, we can get

$$C\boldsymbol{X} = [c_1 \quad \cdots \quad c_q \quad c' \quad \cdots \quad c'] \begin{bmatrix} x_1 \\ \vdots \\ x_q \\ x_{q+1} \\ \vdots \\ x_n \end{bmatrix} = [c_1 \quad \cdots \quad c_q \quad c'] \begin{bmatrix} x_1 \\ \vdots \\ x_q \\ \sum_{i=q+1}^{n} x_i \end{bmatrix} = [c_1 \quad \cdots \quad c_q \quad c'] \begin{bmatrix} x_1 \\ \vdots \\ x_q \\ \tilde{x} \end{bmatrix} = \widetilde{C}\widetilde{\boldsymbol{X}}$$

Similarly, in the constraint of **F2** (i.e., Eq. (9b)), we have

$$G A \boldsymbol{X} = \begin{bmatrix} \sum_{j=1}^{m} g_{1,j} a_{j,1} & \sum_{j=1}^{m} g_{1,j} a_{j,2} & \cdots & \sum_{j=1}^{m} g_{1,j} a_{j,n} \\ \sum_{j=1}^{m} g_{2,j} a_{j,1} & \sum_{j=1}^{m} g_{2,j} a_{j,2} & \cdots & \sum_{j=1}^{m} g_{2,j} a_{j,n} \\ \vdots & \vdots & \ddots & \vdots \\ \sum_{j=1}^{m} g_{K,j} a_{j,1} & \sum_{j=1}^{m} g_{K,j} a_{j,2} & \cdots & \sum_{j=1}^{m} g_{K,j} a_{j,n} \end{bmatrix} \begin{bmatrix} x_1 \\ x_2 \\ \vdots \\ x_n \end{bmatrix}$$

From the second condition, we have $R' = [r'_1, r'_2, \ldots, r'_K]$ such that $\sum_{j=1}^{m} g_{k,j} a_{j,i} = r'_k$ for all $i \in \{q + 1, \ldots, n\}$ and for all $k \in \{1, \ldots, K\}$. Therefore, we can get

$$G A \boldsymbol{X} = \begin{bmatrix} \sum_{j=1}^{m} g_{1,j} a_{j,1} & \cdots & \sum_{j=1}^{m} g_{1,j} a_{j,q} & r'_1 \\ \sum_{j=1}^{m} g_{2,j} a_{j,1} & \cdots & \sum_{j=1}^{m} g_{2,j} a_{j,q} & r'_2 \\ \vdots & \ddots & \vdots & \vdots \\ \sum_{j=1}^{m} g_{K,j} a_{j,1} & \cdots & \sum_{j=1}^{m} g_{K,j} a_{j,q} & r'_K \end{bmatrix} \begin{bmatrix} x_1 \\ \vdots \\ x_q \\ \tilde{x} \end{bmatrix} = \widetilde{A}\widetilde{\boldsymbol{X}}$$

By repeating this process, we can pack commodities into multiple package commodities. As a result, we achieve a new formulation.

**Formulation F3 (Variable Packing)**

Minimize:

$$\mathcal{J} = \tilde{C}\widetilde{\mathcal{X}} \tag{11a}$$

Subject to:

$$\tilde{A}\widetilde{\mathcal{X}} \leq G\boldsymbol{b} \tag{11b}$$

According to the above analysis, the formulation **F3** is equivalent to **F2**.

In summary, we have shown how to find a lower-bound formulation through constraint aggregation and variable packing for general linear programming problems. It is necessary to first find the aggregation matrix $G$ that satisfies the two defining properties (i.e., conditions 1 and 2). Then, we need to identify the variables whose coefficients satisfy the two variable packing conditions (i.e., conditions 3 and 4). This sequence can be generalized to the commodity packing in the space logistics formulation and formulation **F3** can be generalized to the multi-fidelity formulation. Thus, together with the prefixed formulation discussed before, we have the following relationship:

$$\mathcal{J}_{\text{multi\_fidelity}} \leq \mathcal{J}_{\text{full\_size}} \leq \mathcal{J}_{\text{prefixed}}$$

Bounding the computationally prohibitive full-size formulation from both the upper and lower sides enables us to find the approximation solution of the computationally prohibitive full-size formulation with the knowledge about the worst possible approximation error.

**B. Preprocessing Algorithm for Automatic Commodity Packing in Space Logistics**

Although the previous subsection showed an efficient way to pack the commodities in space logistics formulation under certain conditions, we still need a method to identify what commodities are able to be packed in each arc and then find the aggregation matrix to aggregate corresponding constraints so that all conditions are satisfied. Therefore, this subsection proposes a preprocessing algorithm to compile a multi-fidelity optimization formulation automatically for the full-size space infrastructure optimization problem. The consequent formulation performs constraint and variable aggregations in a subset of network arcs, which achieves a lower bound approximation of the original full-size optimization formulation.

Considering a full-size space infrastructure optimization problem as shown in the formulation (1a)-(1d), we can identify the packable commodities leveraging the special structure of this formulation. In the mass balance constraints (i.e., Eq. (1b)), each constraint is designated to guarantee the mass balance of one type of commodity. The commodity transformation matrix $F$ defines the interactions between commodities. To make the commodities packable, they should have the same transformation coefficients with respect to all other commodities. The concurrency constraints (i.e., Eq. (1c)) provide the commodity flow upper bound by considering the total weights of different commodities. For example, the total mass of crew, consumables, instruments, and infrastructure elements must be smaller or equal to the spacecraft total payload capacity; this constraint has a set of packable commodity weights. Therefore, the packable commodities should have the same weight coefficients in all concurrency constraints. After identifying the packable commodities, they can be packed directly in the concurrency constraints without an aggregation matrix. The time window constraints (i.e., Eq. (1d)) are also defined specifically for each type of commodity. However, by definition, the time window is always the same for different commodities in one specific arc. In summary, according to conditions 3 and 4 in Sec. IV.A.2, to pack the packable commodities in space transportation, the associated coefficients must satisfy the following three commodity packing conditions:

1) For the objective function, Eq. (1a), the cost coefficients of packable commodities need to be equal;

2) For the mass balance constraint, Eq. (1b), the transformation coefficients of packable commodities with respect to all other commodities need to be equal;

3) For all concurrency constraints, Eq. (1c), the weight coefficients of packable commodities need to be equal.

Based on the preceding commodity packing conditions, we can propose a preprocessing algorithm to automatically identify the packable commodities and aggregation matrices in the original full-size space infrastructure optimization problem. The pseudo code of the preprocessing is shown as follows. We assume there are $R$ types of commodities in the system. Note that, in this pseudo code, there is a sorting process after identifying packable commodity index sets. The reason for this step is to enable flexible packing decision; if the users prefer to generate fewer package commodities than the number of packable commodity index sets $Q = \{Q_u : u = 1, \dots, U\}$ for an arc (i.e., only $N$ package commodities, where $N \leq U$), they can generate the $N$ most impactful package commodities in the sorted list, where "most impactful" means it contains the most packable commodities. Fewer package commodities, which means fewer commodities are packed, leads to a tighter lower-bound of the optimization.

To generate the aggregation matrix $G$ for the mass balance constraint and time window constraint, we first need to identify the packable commodity index set, denoted by $\check{\zeta}$ as shown in the preprocessing pseudo code. If we assume that we would like to generate $N$ package commodities, then $N = |\check{\zeta}|$. Each subset in $\check{\zeta}$ represents the corresponding packable commodities that will be packed into one package commodity. Suppose that $L$ types of commodities are packed into $N$ package commodities, then $L = \sum_{S_\tau \in \zeta} |S_\tau|$. Therefore, before the commodity packing, the number of variables over each arc is $R$, where each variable represents one type of commodity. After the commodity packing, the number of variables is $K = N + R - L$, where the first $N$ variables represent the package commodities, they contain the information of $L$ types of commodities that are packed; the remaining $R - L$ variables represent commodities that are not packed. Note that the mass balance constraints and the time window constraints are defined for each commodity independently. Therefore, before the commodity packing, the number of mass balance constraints or the time window constraints over each arc is also $R$; after the commodity packing, the number of these constraints becomes $K = N + R - L$.

**Algorithm 1. Preprocessing for commodity packing pseudo code**

---

**For $\forall (v, i, j, t) \in \mathcal{A}$:**

**Step 1. For the cost matrix in the objective function, $c_{vijt}$:** Let $\sigma_1 = \{S_k : k = 1, \ldots, q\}$ be a partition of the commodity indices $\{1, \ldots, R\}$ and define $n_k = |S_k|$. The partition satisfies

$$c_{vijt}[l] = c_{vijt}[l'] \quad \forall l, l' \in S_k \quad \forall k$$

$$\bigcup_{k=1}^{q} S_k = \{1, \ldots, R\} \text{ and } S_k \cap S_{k'} = \emptyset \quad \forall k \neq k'$$

**Step 2. For the transformation matrix in mass balance constraint, $F_{vij}$:** Let $\sigma_2 = \{S_f : f = 1, \ldots, q'\}$ be a partition of the commodity indices $\{1, \ldots, R\}$ and define $n_f = |S_f|$. The partition satisfies

$$F_{vij}[l, u] = F_{vij}[l', u] \quad \forall l, l' \in S_f, \forall u \in \{1, \ldots, R\} \setminus \{l, l'\}, \quad \forall f$$

$$\bigcup_{f=1}^{q'} S_f = \{1, \ldots, R\} \text{ and } S_f \cap S_{f'} = \emptyset \quad \forall f \neq f'$$

**Step 3. For the concurrency matrix in concurrency constraint, $H_{vij}$:** Let $\sigma_3 = \{S_h : h = 1, \ldots, q''\}$ be a partition of the commodity indices $\{1, \ldots, R\}$ and define $n_h = |S_h|$. The partition satisfies

$$H_{vij}[:, l] = H_{vij}[:, l'] \quad \forall l, l' \in S_h \quad \forall h$$

$$\bigcup_{h=1}^{q''} S_h = \{1, \ldots, R\} \text{ and } S_h \cap S_{h'} = \emptyset \quad \forall h \neq h'$$

**Step 4. Find all intersection sets**

$$\zeta = \{S_\tau : \tau = 1, \ldots, U | S_\tau \neq \emptyset \text{ and } S_\tau = S_k \cap S_f \cap S_h, \forall S_k \in \sigma_1, \forall S_f \in \sigma_2, \forall S_h \in \sigma_3\}$$

**Step 5. Identify the packable commodities**

**If $\zeta = \emptyset$:**

   **Step 5.1.** There are no packable commodities in this arc. **Screen the next arc: Go to Step 1.**

---

**Else**:

**Step 5.2.** Define the cardinality $n_\tau = |S_\tau|$.

**Step 5.3.** Perform sorting in a descending order based on the cardinality for $\forall S_\tau \in \zeta$ and get a new set $\tilde{\zeta}$.

**Step 5.4.** Based on the predefined preference, define the number of package commodities as $N$ ($N \leq U$), the packable commodity index set as $\check{\zeta} = \{S_\tau : \tau = 1, \ldots, N | S_\tau \in \tilde{\zeta}\}$.

**Step 6. Find the aggregation matrix for the mass balance constraint and time window constraint:**

**Step 6.1.** Get the number of commodities that will be packed: $L = \sum_{S_\tau \in \check{\zeta}} |S_\tau|$.

**Step 6.2.** Get the number of variables after commodity packing: $K = N + R - L$.

**Step 6.3.** For this arc, define $G = [\boldsymbol{g}_1^T \quad \boldsymbol{g}_2^T \quad \ldots \quad \boldsymbol{g}_K^T]^T$, where each row of the aggregation matrix is a $1 \times R$ weighting vector, $\boldsymbol{g}_k^T$, that satisfies

For $\forall k \in \{1, \ldots, N\}$ (the first $N$ variables are package commodities):

$$\begin{cases} \boldsymbol{g}_k[j] = 1 & if \ j \in S_k \\ \boldsymbol{g}_k[j] = 0 & if \ j \notin S_k \end{cases} \quad S_k \in \check{\zeta}$$

For $\forall k \in \{N+1, \ldots, K\}$ (the remaining variables are for commodities that are not packed):

$$\begin{cases} \boldsymbol{g}_k[j] = 1 & if \ k = N + j - L \\ \boldsymbol{g}_k[j] = 0 & if \ k \neq N + j - L \end{cases}$$

**Step 7. Screen the next arc: Go to Step 1.**

## V. Case Study and Analysis

This section evaluates the performances of the proposed space infrastructure optimization formulations with a case study on a multi-mission human lunar exploration campaign, considering ISRU system designs. The mission scenario, including mission demand, spacecraft design, and ISRU architecture models, is first introduced in Sec. V.A, and then Sec. V.B evaluates the performance of the formulations. Note that although this paper introduces the formulations in the order of the prefixed, full-size, multi-fidelity formulations, the later analysis considers the full-size optimization formulation as the baseline and compares the other two formulations against it; this is because the full-size optimization is the most accurate and computationally costly one, and we are interested in the solution quality and the computational cost of the prefixed formulation (i.e., the upper-bound formulation) and the multi-fidelity formulation (i.e., the lower-bound formulation).

### A. Mission Scenario

A simple scenario is considered as a case study where all formulations (including the full-size formulation) can complete its computation within a reasonable time. We consider a cis-lunar transportation system with Earth, low-Earth orbit (LEO), geosynchronous equatorial orbit (GEO), Earth-Moon Lagrangian point 1 (EML1), and the Moon.

The five-node transportation network model and the spaceflight $\Delta V$ values are shown in Fig. 3. Note that we do not consider the propellant cost from Earth to LEO (i.e., the $\Delta V$ is considered as zero); instead, Earth is assumed as the main supply node and the arc from Earth to LEO is convenient to calculate the space mission cost. Every year, 5 astronauts fly to the Moon with habitat and equipment. These demands are considered as one type of general payload together with a crew cabin. The total mass of the crew cabin and lunar equipment is assumed as 30,000 kg, which is estimated based on the Apollo mission [25]. The astronauts stay on the lunar surface for 120 days and then come back with lunar samples and materials. The total mass of crew cabin and lunar samples is assumed as 5,000 kg and they are delivered back to the Earth at the end of the mission. For this mission design, the optimizer needs to decide whether it requires ISRU systems to support the transportation, whether the system needs a propellant depot, and where we should deploy the depot (i.e., LEO, GEO, or EML1) if needed. We assume that a spacecraft can serve as a propellant depot if it stays at a node during the mission [12]. The mission demands and supplies are summarized in Table 1. Note that, the mission demands and supplies are defined at the same time step for each flight to minimize the number of time steps assigned for the transportation.

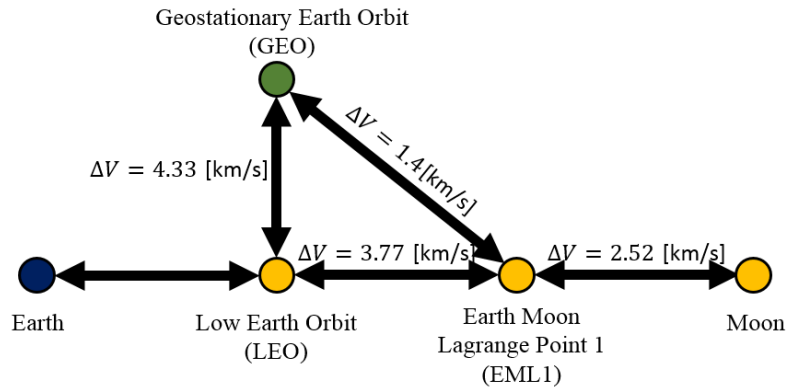

**Fig. 3 Cis-lunar transportation network model.**

**Table 1 Lunar exploration demands and supplies**

| Payload Type | Node | Time, day | Supply |
|---|---|---|---|
| Go to the Moon | | | |
| ISRU, propellant & food, kg | Earth | All the time | +∞ |
| Crew cabin & equipment, kg | Earth | 240 | +30,000 [25] |
| Crew cabin & equipment, kg | Moon | 240 | -30,000 [25] |
| Back to Earth | | | |
| Crew cabin & lunar sample, kg | Moon | 360 | +5,000 |
| Crew cabin & lunar sample, kg | Earth | 360 | -5,000 |

We need spacecraft to deliver payloads from Earth to the Moon. To simplify the analysis, the spacecraft design is not considered as part of the trade space in space logistics optimization. Instead, two types of spacecraft with fixed design parameters are considered for space transportation. Spacecraft 1 is modeled based on the Advanced Cryogenic Evolved Stage (ACES) from United Launch Alliance [12]. It uses liquid hydrogen and liquid oxygen (LH2/LOX) as the propellant. The spacecraft structure mass is 5,917 kg and the propellant tank capacity is 68,040 kg [12]. Because of the implementation of long-duration storage technologies in ACES propellant tanks, the LH2/LOX propellant boiloff rate is considered as zero during space transportation. Spacecraft 2 is modeled based on the lunar surface access module (LSAM) descent stage pressure-fed design from the green propellants study. The design parameters are found in the SpaceNet database [26]. It uses liquid methane and liquid oxygen (LCH4/LOX) as the propellant. The spacecraft design assumptions are listed in Table 2. For simplicity, we assume that both spacecraft can be used for all trajectories in the transportation network, including lunar landing and ascending. Also, they are considered as single-stage transportation vehicles, but they can be combined to form a larger transportation vehicle.

**Table 2 Spacecraft design parameters.**

| Parameter | Assumed value |
|---|---|
| **Spacecraft 1** | |
| Propellant type | LH2/LOX |
| Propellant capacity, kg | 68,040 [12] |
| Structure mass, kg | 5,917 [12] |
| Propellant $I_{sp}$, s | 420 |
| Propellant component mass ratio | $O_2:H_2$=5.5:1 |
| **Spacecraft 2** | |
| Propellant type | LCH4/LOX |
| Propellant capacity, kg | 40,737 [26] |
| Structure mass, kg | 6,560 [26] |
| Propellant $I_{sp}$, s | 350 |
| Propellant component mass ratio | $O_2:CH_4$=3.5:1 |

The ISRU infrastructure design model is another essential part of the space infrastructure optimization case study. For the lunar exploration campaign considered in this paper, the ISRU architecture design models are listed in Table 3. These models are extrapolated from historical ISRU infrastructure design concept literature and prototypes by Chen et al. [27].

In Table 3, the reference product is used to size the ISRU subsystems. For reactors and excavators, the specific power and specific mass mean the power demand and the system mass needed to reach 1 kg/hr productivity of the reference product. For storage systems and power systems, the specific power and specific mass mean the necessary system size to store 1 kg resource, 1 kWh energy or to supply 1 kW power. The soil/water extraction process and the excavator are classified based on different soil types, soil @3% $H_2O$ and soil @8% $H_2O$. Note that the regolith water concentration values assumed here are extrapolated from the literature relating to Martian surface soil [27]. They are used as example values only here. Because of the difference in lunar regolith composition, the hydrogen reduction process also has different productivity in different regions. Moreover, according to the ISRU infrastructure design prototype, we assume that rigid solar concentrators provide thermal energy to the HR and CR reactors [4]. They are considered as part of the reactors. Therefore, the nominal power demands of the HR and CR reactors are zero. In this case, we only consider ISRU systems for $O_2$ and $H_2$ generation during the mission. It is up to the optimizer's choice whether to use Spacecraft 1 and leverage ISRU systems or to use Spacecraft 2 and deliver all necessary propellant from Earth.

Besides the ISRU infrastructure sizing models, mission operation management is also critical to be considered in space logistics optimization. It includes rocket launch frequency, ISRU system maintenance [14-15], power system working environment, degradation, and energy storage efficiencies [28-32]. The mission operation assumptions are listed in Table 4. The rocket launch interval determines the frequency of mission operation. We define that the mission operation time windows are open for a few time steps after each rocket launch opportunity. When the mission operation time windows are closed, space flights are not permitted. The ISRU maintenance rate means that every year, the mass of maintenance spare demand is equivalent to 10% of the ISRU system total mass [14-15].

Table 3 ISRU infrastructure design models. [27]

| System | Chemistry reactions | Reference product | Specific power, kW | Specific mass, kg |
|---|---|---|---|---|
| **Reactor** | | | | |
| Soil/Water extraction (SWE) | $Soil \rightarrow H_2O$ | $H_2O$, kg/hr | @3%: 13.7 @8%: 7 | @3%: 357 @8%: 195 |
| Direct water electrolysis (DWE) | $2H_2O \rightarrow 2H_2 + O_2$ | $O_2$, kg/hr | 5.83 | 83.3 |
| Molten regolith electrolysis (MRE) | $Soil \rightarrow O_2$ | $O_2$, kg/hr | 26.94 | 197.58 |
| Hydrogen reduction (HR) | $Soil + H_2 \rightarrow H_2O$ | $H_2O$, kg/hr | 0 | @equator: 228 @pole: 482 |
| Carbothermal reduction (CR) | $Soil + 2CH_4 + 2H_2 \rightarrow 2CH_4 + 2H_2O$ | $H_2O$, kg/hr | 0 | 520.5 |
| **Soil extraction system** | | | | |
| Excavator for soil @3% $H_2O$ | —— | Soil, kg/hr | 0.004 | 0.38 |
| Excavator for soil @8% $H_2O$ | —— | Soil, kg/hr | 0.027 | 23 |
| **Storage system** | | | | |
| $O_2$ storage | —— | $O_2$, kg | 0.0088 | 5.15 |
| $H_2$ storage | —— | $H_2$, kg | 0.0267 | 3.33 |
| $H_2O$ storage | —— | $H_2O$, kg | 0 | 40 |
| $CH_4$ storage | —— | $CH_4$, kg | 0.0073 | 1.67 |
| **Power system** | | | | |
| Photovoltaic (PV) power system | —— | Power, kW | —— | 6.8 (@ 1 AU) |

| Energy storage system: battery | —— | Energy, kWh | —— | 4 |
| Energy storage system: fuel cell | —— | Energy, kWh | —— | 2 |
| Fission surface power system (FSPS) | —— | Power, kW | —— | 150 |
| Radioisotope power system (RPS) | —— | Power, kW | —— | 124 |

**Table 4 Mission operation parameters and assumptions.**

| Parameter | Assumed value |
|---|---|
| Rocket launch interval, day | 120 |
| ISRU maintenance, system mass/yr | 10% [14-15] |
| Solar irradiance (@ 1 AU), kW/m$^2$ | 1.36 [28] |
| PV radiation degradation, /sol | 0.014% [29] |
| Battery charging efficiency | 95% [30] |
| Fuel cell energy efficiency | 60% [31] |
| RPS degradation rate, /yr | 1.9% [32] |

The problem is solved using the Gurobi 8.1 solver through Python on an i9-9900k, 3.6GHz platform with 32GB RAM. The detailed analysis and discussion of this human lunar exploration campaign case study are shown in the next section.

### B. Comparison of Optimization Formulations

This section compares the solution and computational cost of the prefixed infrastructure optimization formulation (i.e., the upper-bound formulation), the full-size infrastructure optimization formulation (i.e., the baseline formulation), and the proposed multi-fidelity optimization formulation (i.e., the lower-bound formulation). We consider a lunar exploration campaign with multiple consecutive lunar missions, with a mission operation frequency of 120 days. The lunar landing area is in the equatorial region with lunar regolith @3% $H_2O$. The initial mass in low-Earth orbit (IMLEO) is used as the mission cost metric. It is a widely used mission cost measurement in past space logistics optimization literature [13-15]. As a baseline mission scenario, the FSPS is selected as the stationary power supply system on the lunar surface. The PV power system and energy storage system are considered as candidate power sources in space.

By fixing the number of human lunar missions to three and changing other mission scenario parameters, we can evaluate the performance of three optimization formulations under different settings. The ISRU infrastructure design models shown in Table 3 are relatively conservative models. With the development of technology and material science, ISRU systems can have higher productivity and lower system structure mass. Table 5 compares the infrastructure optimization formulation performances with respect to ISRU productivity. It shows the results when the ISRU productivities are 100%, 125%, and 150% of the original design models. The mission cost errors illustrate the mission cost difference of solutions with respect to the results of the baseline full-size optimization formulation.

In Table 5, we can find that multi-fidelity optimization can provide a very accurate approximation of the full-size formulation at a significant computational cost reduction. The computation time reduction was more than 60%, whereas the performance loss is within 2.5%. This is enabled by packing commodity variables and eliminating infrastructure subsystem tradeoffs during space flights. The (small) solution difference between the multi-fidelity and full-size formulations is caused by the inability of the multi-fidelity formulation to distinguish different commodity types when they are packed together; for example, when two commodities are packed and then unpacked later on, we lose the information about the original mass ratio between these two commodities, which can lead to an overoptimistic solution.

On the other hand, the upper-bound solutions provided by the prefixed optimization formulation is much larger than optimal solutions. The physical meaning of the prefixed infrastructure optimization is that it ignores the infrastructure subsystem trade studies and their interactions with space mission planning. It considers the infrastructure as an integrated system. We can still size the infrastructure; however, the mass ratios between infrastructure subsystems are fixed in advance before considering space logistics. Therefore, it can provide an upper-bound, feasible solution, which is significantly larger than the optimal solution. It is also the fastest method among three infrastructure optimization formulations because it has the least variables and constraints and explores the smallest design space.

**Table 5 Comparison of formulation performances with respect to ISRU productivity.**

| ISRU productivity index | Optimization formulation | Mission cost (IMLEO), kg | Mission cost errors | Computation time, s | Computation time reduction |
|---|---|---|---|---|---|
| 100% (default) | Prefixed (Upper Bound) | 565,622.9 | 33.7% | 110.9 | -89.6% |
| | Full-size (Baseline) | 422,930.7 | —— | 1,062.3 | —— |
| | Multi-fidelity (Lower Bound) | 414,393.7 | -2.0% | 92.5 | -91.3% |
| 125% | Prefixed (Upper Bound) | 528,563.7 | 33.6% | 91.1 | -92.0% |
| | Full-size (Baseline) | 395,422.6 | —— | 1,135.7 | —— |
| | Multi-fidelity (Lower Bound) | 394,302.7 | -0.3% | 301.9 | -73.4% |
| 150% | Prefixed (Upper Bound) | 513,612.6 | 37.1% | 25.7 | -95.7% |

| | Full-size (Baseline) | 374,732.9 | —— | 592.2 | —— |
| | Multi-fidelity (Lower Bound) | 366,229.3 | -2.3% | 246.8 | -58.3% |

We can vary the problem complexity by changing the number of human lunar missions or the rocket launch frequency. If we fix the ISRU productivity as normal and increase the number of human lunar missions from 3 to 4 and 5, the mission planning results are shown in Table 6. It shows that the mission cost errors of the multi-fidelity optimization with respect to the full-size optimization are within 2%. Both the multi-fidelity optimization and the prefixed optimization formulations are significantly faster than the full-size optimization formulation (i.e., >90% computation time reduction).

**Table 6 Optimization formulation performance comparison.**

| Number of human lunar missions | Optimization formulation | Mission cost (IMLEO), kg | Mission cost errors | Computation time, s | Computation time reduction |
|---|---|---|---|---|---|
| 3 (default) | Prefixed (Upper Bound) | 565,622.9 | 33.7% | 110.9 | -89.6% |
| | Full-size (Baseline) | 422,930.7 | —— | 1,062.3 | —— |
| | Multi-fidelity (Lower Bound) | 414,393.7 | -2.0% | 92.5 | -91.3% |
| 4 | Prefixed (Upper Bound) | 671,716.3 | 31.3% | 321.3 | -96.7% |
| | Full-size (Baseline) | 511,476.6 | —— | 9,607.4 | —— |
| | Multi-fidelity (Lower Bound) | 509,792.9 | -0.3% | 439.8 | -95.4% |
| 5 | Prefixed (Upper Bound) | 774,626.1 | 29.7% | 693.8 | -98.4% |
| | Full-size (Baseline) | 597,300.8 | —— | 42,675.8 | —— |
| | Multi-fidelity (Lower Bound) | 596,347.8 | -0.2% | 2,074.1 | -95.1% |

We can also fix the number of human lunar missions to three and the ISRU productivity as normal, then change the launch frequency to evaluate its impact on ISRU infrastructure design, especially the storage system design. As there is a 120-day long human lunar exploration at the end of each year, the human lunar mission begins on day 240 in each year. By varying the rocket launch frequency interval to 60, 120 (default), or 240 days, there are 3, 1 or 0 extra cargo mission opportunities before each human lunar mission. The formulation performance comparison under different launch frequencies is shown in Table 7.

Table 7 shows that the performance of the multi-fidelity optimization formulation is stable. The mission cost errors are always within 2% compared with the optimal solutions from the full-size optimization formulation. If we observe the computation times in Table 6 and Table 7, we can analyze a general trend in the computational time saving by the multi-fidelity formulation. In Table 6, as the number of human lunar missions increases, the computation time reduction of the multi-fidelity formulation increases slightly from 91% to 95%. In Table 7, as the rocket launch opportunity interval decreases (i.e., from 240 to 60), the time steps considered in the optimization increase significantly, and the computation time reduction of the multi-fidelity formulation increases from 70% to more than 95%. These observations show that the proposed multi-fidelity optimization formulation achieves a large computational time saving compared with the full-size formulation for complex space mission design problems.

**Table 7 Comparison of formulation performances with respect to the launch frequency.**

| Launch frequency, day | Optimization formulation | Mission cost (IMLEO), kg | Mission cost errors | Computation time, s | Computation time reduction |
|---|---|---|---|---|---|
| 240 | Prefixed (Upper Bound) | 697,800.9 | 65.0% | 7.1 | -94.8% |
| | Full-size (Baseline) | 422,930.7 | —— | 135.9 | —— |
| | Multi-fidelity (Lower Bound) | 414,393.7 | -2.0% | 39.9 | -70.6% |
| 120 (default) | Prefixed (Upper Bound) | 565,622.9 | 33.7% | 110.9 | -89.6% |
| | Full-size (Baseline) | 422,930.7 | —— | 1,062.3 | —— |
| | Multi-fidelity (Lower Bound) | 414,393.7 | -2.0% | 92.5 | -91.3% |
| 60 | Prefixed (Upper Bound) | 480,705.3 | 13.7% | 604.9 | -98.2% |
| | Full-size (Baseline) | 422,926.5 | —— | 33,383.0 | —— |
| | Multi-fidelity (Lower Bound) | 414,388.1 | -2.0% | 1,581.9 | -95.3% |

Moreover, the results also show that the launch frequency and the sizing of infrastructure storage systems need to be considered concurrently to find the optimal infrastructure design. With a higher launch frequency, a smaller storage system is needed because resources produced by the infrastructure can be delivered to other destinations through spacecraft when mission time windows are open. Keeping this intuition in mind, the storage system design in the prefixed optimization formulation is pre-set to be able to store the exact amount of resources produced between two mission time windows. For example, if the launch frequency is 240 days, then the storage system in the prefixed optimization formulation is exactly able to store the resources produced in 240 days. In Table 7, the mission cost results by the full-size and the multi-fidelity optimization formulations show that the launch frequency has limited influence on mission costs for this mission scenario. However, the mission cost from the prefixed optimization formulation decreases significantly as the launch frequency increases, which leads to a decrease in infrastructure storage system size. This result shows that this mission scenario may prefer small infrastructure storage systems.

To confirm this hypothesis, we conduct a sensitivity analysis on the ISRU storage system sizing under the default launch frequency (i.e., 120 days). The results are shown in Table 8. We find that as we decrease the storage system size, the mission costs obtained through the prefixed optimization formulation decrease dramatically until the storage system is too small to make the mission feasible. Note that our full-size formulation's solution is still much better than any of the prefixed formulations tested here. This result shows that our proposed interdisciplinary space infrastructure optimization methods can optimize the ISRU storage size as well as any other ISRU subsystems by concurrently capturing the detailed interactions between each infrastructure subsystem and space transportation mission planning in an optimal way. The optimal subsystem designs cannot be achieved by considering space infrastructure design independently in advance and treating it as a black box in space logistics.

**Table 8 Sensitivity analysis of ISRU storage system sizing.**

| Optimization formulation | ISRU storage system size | Mission cost (IMLEO), kg | Mission cost errors | Computation time, s | Computation time reduction |
|---|---|---|---|---|---|
| Prefixed (Upper Bound) | 100% | 565,622.9 | 33.7% | 110.9 | -89.6% |
|  | 80% | 524,412.9 | 24.0% | 41.5 | -96.1% |
|  | 60% | 494,243.2 | 16.9% | 37.2 | -96.5% |
|  | 40% | 467,241.3 | 10.5% | 59.4 | -94.4% |
|  | 20% | 444,414.1 | 5.1% | 151.2 | -85.8% |
|  | 0% | *infeasible* | —— | —— | —— |
| Full-size (Baseline) | —— | 422,930.7 | —— | 1,062.3 | —— |
| Multi-fidelity (Lower Bound) | —— | 414,393.7 | -2.0% | 92.5 | -91.3% |

## VI. Conclusion

This paper proposes a system-level space infrastructure and logistics mission design optimization framework to perform architecture trade studies. A new space infrastructure logistics optimization problem formulation is proposed that considers infrastructure subsystems' internal interactions and their external synergistic effects with space logistics simultaneously. A natural implementation of this formulation is referred to as the full-size formulation, which explores a larger trade space and thus provides the same or a better (i.e., lower-cost) solution than the traditional prefixed formulation. However, the inherent limitation of this full-size formulation is its prohibitive computational cost for complex systems. In response to this challenge, another new multi-fidelity optimization formulation is developed by varying the granularity of the commodity type definition over the network graph. The developed multi-fidelity formulation can find an approximation lower-bound solution to the full-size problem computationally efficiently with little sacrifice in the solution quality. A multi-mission human lunar exploration campaign case study shows the consistent improvement of the multi-fidelity optimization formulation in computational efficiency. For the tested cases, the multi-fidelity optimization formulation found solutions that are within 2-3% of those of the full-size optimization formulation with a significant computational time reduction (>90% for the majority of the tested cases). The sensitivity analysis of launch frequency demonstrates the value of the proposed interdisciplinary infrastructure optimization method.

Future research can include the consideration of uncertainties in space mission planning to evaluate the performance and technology reliability under stochastic mission scenarios. Further implementations can also be explored to consider technology trade studies for life support systems or scientific instruments in space logistics optimization.

## Funding Sources

This material is partially based upon work supported by the funding from the NASA NextSTEP program (80NSSC18P3418) awarded to the University of Illinois, where the original version of this work was initiated. Any opinions, findings, and conclusions or recommendations expressed in this material are those of the authors and do not necessarily reflect the views of the National Aeronautics and Space Administration.